\documentclass{amsart}

\usepackage{a4wide, amssymb, mathrsfs,stmaryrd}
\usepackage{fancyhdr}
\usepackage[T1]{fontenc}

\usepackage[arrow,matrix,curve]{xy}

\def\id{\textrm{id}}

\def\xto#1{\xrightarrow{#1}}

\newcommand{\nc}{\newcommand}

\nc{\sA}{{\mathscr A}}
\nc{\wt}{\widetilde} \nc{\bl}{\bullet} \nc{\al}{\alpha}
\nc{\sg}{\sigma} \nc{\vf}{\varphi} \nc{\om}{\omega}
\nc{\ve}{\varepsilon} \nc{\ol}{\overline} \nc{\lb}{\lambda}
\nc{\Lb}{\Lambda} \nc{\Gm}{\Gamma} \nc{\cP}{{\mathscr P}}
\nc{\sB}{{\mathscr B}}
\nc{\ul}{\underline} \nc{\os}{\overset} \nc{\us}{\underset}
\nc{\pa}{\partial} \nc{\wh}{\widehat} \nc{\sbs}{\subset}
\nc{\lra}{\longrightarrow} \nc{\all}{\allowdisplaybreaks}
\nc{\Ker}{\operatorname{Ker}} \nc{\Img}{\operatorname{Im}}
\nc{\Kan}{\operatorname{Kan}} \nc{\Hom}{\operatorname{Hom}}
\nc{\Imm}{\operatorname{Im}}
\nc{\idd}{\operatorname{id}}
\nc{\Top}{\operatorname{\textbf{Top}}}
\nc{\Ab}{\operatorname{\textbf{Ab}}}

\newtheorem{theo}{Theorem}[section]
\newtheorem{prop}[theo]{Proposition}
\newtheorem{lem}[theo]{Lemma}
\newtheorem{coro}[theo]{Corollary}

\theoremstyle{definition}
\newtheorem{defi}[theo]{Definition}
\newtheorem{exmp}[theo]{Example}
\newtheorem{remk}[theo]{Remark}

\newcommand{\mybox}{\ensuremath \Box}

\newenvironment{prf}{\noindent\textbf{Proof.}}{\mybox}

\usepackage{footmisc}
\DefineFNsymbols{void}{{\hspace{-1cm}}}
\setfnsymbol{void}

\begin{document}

\def\L{\ensuremath{{\mathbf{L}}}}
\def\Z{{\mathbb Z}}

\def\D{{\mathscr D}}
\def\B{{\mathscr B}}
\def\Q{{\mathscr Q}}
\def\M{{\mathscr M}}
\def\D{{\mathscr D}}
\def\R{{\mathscr R}}
\def\E{{\mathscr E}}
\def\K{{\mathscr K}}
\def\W{{\mathscr W}}
\def\N{{\mathscr N}}
\def\T{{\mathscr T}}
\def\A{{\mathscr A}}
\def\C{{\mathscr C}}
\def\P{{\mathscr P}}
\def\V{{\mathscr V}}
\def\G{{\mathscr G}}
\let\xto\xrightarrow

\advance\baselineskip by.3\baselineskip

\title[]
{COMPARISON OF CUBICAL AND SIMPLICIAL DERIVED FUNCTORS}

\author{Irakli Patchkoria}

\begin{abstract}
In this note we prove that the simplicial derived functors introduced by Tierney and Vogel [TV69] are naturally isomorphic to the cubical derived functors introduced by the author in [P09]. We also explain how this result generalizes the well-known fact that the simplicial and cubical singular homologies of a topological space are naturally isomorphic.
\end{abstract}
\maketitle
\setcounter{section}{0}
\section{Introduction}

In [TV69] Tierney and Vogel for any functor $F \colon \C \lra \B$, where $\C$ is a category with finite limits and a projective class $\P$, and $\B$ is an abelian category, constructed simplicial derived functors and investigated relationships of their theory with other theories of derived functors. Namely, they showed that if $\C$ is abelian and $F$ is additive, then their theory coincides with the classical relative theory of Eilenberg-Moore [EM65], whereas if $\C$ is abelian and $F$ is an arbitrary functor, then it gives a generalization of the theory of Dold-Puppe [DP61]. Besides, they proved that their derived functors are naturally isomorphic to the cotriple derived functors of Barr-Beck ([BB66], [BB69]) if there is a cotriple in $\C$ that realizes the given projective class $\P$.

The key point in the construction of the derived functors by
Tierney and Vogel is that using $\P$-projective objects and
simplicial kernels, for every $C$ from $\C$ a $\P$-projective
pseudosimplicial resolution can be constructed, which is a $C$-augmented
pseudosimplicial object in $\C$ and which for a given $C$ is
unique up to a presimplicial homotopy.

In [P09] using pseudocubical resolutions instead of pseudosimplicial ones we constructed cubical derived functors for any functor $F \colon \C \lra \B$, where $\C$ is a category with finite limits and a projective class $\P$, and $\B$ is an abelian category. It was shown that if $\C$ is an abelian category, $F$ an additive functor, and $\P$ is closed, then our theory coincides with the theory of Eilenberg-Moore [P09, 4.4]. However, there remained an open question whether the Tierney-Vogel simplicial derived functors and our cubical derived functors are isomorphic in general or not. In this paper we give a positive answer to this question. More precisely, we prove the following

\begin{theo} \label{mtavari} Suppose $\C$ is a category with finite limits, $\P$ a projective class in $ \C$ in the sense of [TV69, \S 2], $\B$ an abelian category, and $F \colon \C \lra \B$ a functor. Let $\L^{\Delta}_nF \colon \C \lra \B$, $n \geq 0$, be the Tierney-Vogel simplicial derived functors of $F$, and $\L^{\Box}_nF \colon \C \lra \B$, $n \geq 0$, the cubical derived functors of $F$. Then there is an isomorphism
$$\xymatrix{\L^{\Delta}_nF(C) \cong \L^{\Box}_nF(C),\;\;\; C \in \C, \;\; n \geq 0,}$$
which is natural in $F$ and in $C$. \end{theo}

The main idea of the proof goes back to Barr and Beck [BB69]. The point is that passing to the unique additive extension 
$$ F_{ad} \colon \Z\C \lra \B$$
of the functor $F$, where $\Z\C$ denotes the free preadditive category generated by $\C$, one verifies that the Eilenberg-Moore derived functors of $F_{ad}$ (with respect to the class $\P$) restricted to $\C$ are naturally isomorphic to the simplicial derived functors of $F$ on the one hand and to the cubical derived functors of $F$ on the other hand.

The paper is organized as follows. In Section 2 the relative Eilenberg-Moore derived functor theory of additive functors is reviewed from [EM65]. In Section 3 we recall the theory of Tierney-Vogel and prove that the simplicial derived functors of $F \colon \C \lra \B$ are just the Eilenberg-Moore derived functors of $F_{ad} \colon \Z\C \lra \B$ restricted to $\C$. Section 4 is devoted to the definition and properties of pseudocubical normalization functor for an idempotent complete preadditive category. Note that the pseudocubical normalization is the main technical tool used in Section 5 to prove that the cubical derived functors of $F \colon \C \lra \B$ are naturally isomorphic to the Eilenberg-Moore derived functors of $F_{ad} \colon \Z\C \lra \B$ restricted to $\C$. In the final section we briefly indicate that Theorem \ref{mtavari} generalizes the classical fact that the simplicial and cubical singular homologies of a topological space are naturally isomorphic.

\vspace{0.2cm}

\section{Partially defined Eilenberg-Moore derived functors}

The following definitions are well-known.

\begin{defi}\label{pread} A preadditive category is a category $\A$ together with the following data:

{\rm (i)} For any objects $X,Y$ in $\A$, the set of morphisms $\Hom_{\A}(X,Y)$ is an abelian group;

{\rm (ii)} For any morphisms $f,g \colon X \longrightarrow Y$, $h \colon W \longrightarrow X$ and $u \colon Y \longrightarrow Z$ in $\A$, the following hold

$$(f+g)h=fh+gh, \;\;\; u(f+g)=uf+ug.$$
\end{defi}

In other words, a preadditive category is just a ring with several objects in the sense of [M72]. 

\begin{defi} \label{complex} Let $\A$ be a preadditive category. An augmented chain complex over an object $A \in \A$ (or just a complex over $A$) is a sequence
$$\xymatrix{ \cdots \ar[r] & C_n \ar[r]^{\partial_n} & C_{n-1} \ar[r] & \cdots \ar[r] & C_2 \ar[r]^{\partial_2} & C_1 \ar[r]^{\partial_1} & C_0 \ar[r]^{\partial_0} & A }$$
such that $\partial_n \partial_{n+1}=0$, $n \geq 0$. 
\end{defi}

\begin{defi} Let $\A$ be a preadditive category and $\P$ a class of objects in $\A$ (which need not be a ``projective class'' in any sense). A complex
$$\xymatrix{\cdots \ar[r] & C_2 \ar[r]^{\partial_2} & C_1 \ar[r]^{\partial_1} & C_0 \ar[r]^{\partial_0} & A }$$
over $A \in \A$ is said to be $\P$-acyclic if for any $Q \in \P$ the sequence of abelian groups 
$$\xymatrix{\cdots \ar[r] & \Hom_{\A}(Q, C_1) \ar[r]^-{{\partial_1}_{*}} & \Hom_{\A}(Q,C_0) \ar[r]^-{{\partial_0}_{*}} & \Hom_{\A}(Q, A) \ar[r] & 0}$$
is exact. \end{defi}

\begin{defi}\label{res} Let $\A$ be a preadditive category and $\P$ a class of objects in $\A$. A $\P$-resolution of an object $A \in \A$ is a $\P$-acyclic complex 
$$\xymatrix{ \cdots \ar[r] & P_2 \ar[r]^{\partial_2} & P_1 \ar[r]^{\partial_1} & P_0 \ar[r]^{\partial_0} & A }$$ 
over $A$ with $P_n \in \P$, $n \geq 0$. \end{defi}

Note that an object $A \in \A$ need not necessarily possess a $\P$-resolution.

There is a comparison theorem for $\P$-resolutions which can be proved using the standard homological algebra arguments (see e.g. [W94, 2.2.7]). More precisely, the following is valid.

\begin{prop}[Comparison theorem]\label{comparison} Let $P_* \longrightarrow A$ be a complex over $A \in \A$ consisting of objects of $\P$, and let $S_* \longrightarrow B$ be a $\P$-acyclic complex. Then any morphism $f \colon A \longrightarrow B$ can be extended to a morphism of augmented chain complexes
$$\xymatrix{ P_* \ar[r] \ar[d]_{\overline{f}} &  A \ar[d]^f \\ S_* \ar[r] & B.}$$
Moreover, any two such extensions are chain homotopic.
\end{prop}

Suppose $\A$ is a preadditive category, $\P$ a class of objects in $\A$, $\B$ an abelian category, $F \colon \A \longrightarrow \B$ an additive functor, and $\A'$ the full subcategory of those objects in $\A$ which possess $\P$-resolutions. Recall that Proposition \ref{comparison} allows one to construct the left derived functors $\L_n^{\P}F \colon \A' \longrightarrow \B$, $n \geq 0$, of $F$ with respect to the class $\P$ as follows. If $A \in \A'$, choose (once and for all) a $\P$-resolution $P_* \longrightarrow A$ and define 
$$\L_n^{\P}F(A)= H_n(F(P_*)), \;\;\;\; n \geq 0.$$

\begin{remk} If $\P$ is a projective class in the sense of [EM65], then $\L_n^{\P}F$, $n \geq 0$, are exactly the derived functors introduced in [EM65, I.3]. Note that in this case $\A'=\A$, i.e., the functors $\L_n^{\P}F$ are defined everywhere. \end{remk}

Further we recall

\begin{defi}[{[EM65, I.2]}] Let $\A$ be a preadditive category and $\P$ a class of objects of $\A$. A sequence
$$\xymatrix{ X \ar[r]^f & Y \ar[r]^g & Z}$$
in $\A$ is said to be $\P$-exact if $gf=0$ and the sequence of abelian groups
$$\xymatrix{ \Hom_{\A}(P,X) \ar[r]^-{f_*} & \Hom_{\A}(P,Y) \ar[r]^-{g_*} & \Hom_{\A}(P,Z)}$$
is exact for any $P \in \P$. \end{defi}

\begin{defi}[{[EM65, I.2]}] A closure of a class $\P$, denoted by $\overline{\P}$, is the class of all those objects $Q \in \A$ for which
$$\xymatrix{ \Hom_{\A}(Q,X) \ar[r]^-{f_*} & \Hom_{\A}(Q,Y) \ar[r]^-{g_*} & \Hom_{\A}(Q,Z)}$$
is exact whenever $\xymatrix{ X \ar[r]^f & Y \ar[r]^g & Z}$ is $\P$-exact. \end{defi}

Clearly, $\P \subseteq \overline{\P}$ and $\P$-exactness is equivalent to $\overline{\P}$-exactness. In particular, $\overline{\overline{\P}}=\overline{\P}$.

Note that if a preadditive category $\A$ has a terminal object, then any $\P$-resolution is a $\overline{\P}$-resolution as well. This together with \ref{comparison} implies the following

\begin{prop} \label{closure} Let $\A$ be a preadditive category with a terminal object, $\P$ a class of objects in $\A$, $\B$ an abelian category, $F \colon \A \longrightarrow \B$ an additive functor, and $A$ an object in $\A$ which possesses a $\P$-resolution. Then there is a natural isomorphism
$$\L_n^{\P}F(A) \cong \L_n^{\overline{\P}}F(A), \;\;\; n \geq 0.$$ \end{prop}

\vspace{0.2cm}

\section{Simplicial derived functors and Eilenberg-Moore derived functors}

In this section we briefly review the construction of simplicial derived functors from [TV69, \S 2] and show that they can be obtained as derived functors of an additive functor.

Let us recall the following definitions.

\begin{defi} A presimplicial object $S$ in a category $\C$ is a family of objects $(S_n \in \C)_{n \geq 0}$ together with morphisms
$$\partial_i \colon S_n \longrightarrow S_{n-1}, \;\;\; n \geq 1, \;\; 0 \leq i \leq n,$$
in $\C$ satisfying the presimplicial identities $$\partial_i \partial_j = \partial_{j-1} \partial_i, \;\;\; i < j.$$ \end{defi}

\begin{defi} Let $S$ be a presimplicial object in a preadditive category $\A$. The unnormalized chain complex $K(S)$ associated to $S$ is defined by 
$$ K(S)_n = S_n,  \;\;\; n \geq 0,$$
$$\partial= \sum_{i=0}^n(-1)^{i} \partial_i \colon K(S)_n \longrightarrow K(S)_{n-1}, \;\;\; n >0.$$ \end{defi}

\noindent The presimplicial identities imply that $\partial^2 =0$.

Now let $\C$ be a category with finite limits, $\P$ a projective class in $\C$ in the sense of [TV69, \S 2], $\B$ an abelian category, and $F \colon \C \longrightarrow \B$ a functor. The simplicial derived functors $\L^{\Delta}_n F$ of $F$ with respect to the class $\P$ are defined as follows. For any object $C \in \C$, choose (once and for all) a $\P$-projective presimplicial resolution 
$$ S \longrightarrow C$$
of $C$ (i.e., a $\P$-exact presimplicial object $S$ augmented over $C$ with $S_n \in \P, \;\; n \geq 0$) and define 
$$\L^{\Delta}_nF(C) = H_n(K(F(S))), \;\;\; n \geq 0.$$
\noindent By the comparison theorem for projective presimplicial resolutions [TV69, (2.4) Theorem], the objects $\L^{\Delta}_nF(C)$ are well-defined and functorial in $F$ and $C$. 

We will now show that the derived functors $\L^{\Delta}_nF$ can be obtained as derived functors of some additive functor. First recall

\begin{lem} \label{contract} Let $S \longrightarrow S_{-1}$ be an augmented presimplicial set. Suppose that $\partial_0 : S_0 \longrightarrow S_{-1}$ is surjective and the following extension condition holds: For any $n \geq 0$ and any collection of $n +2$ elements $x_i \in S_n$, $0 \leq i \leq n+1$, satisfying 
$$\partial_i x_j = \partial_{j-1} x_i, \;\; 0 \leq i < j \leq n+1,$$ 
there exists $x \in S_{n+1}$ such that $\partial_ix=x_i$, $0 \leq i \leq n+1$. Then the augmented chain complex
$$\xymatrix{ K(\Z[S]) \ar[r]^{\partial_0} & \Z[S_{-1}]}$$
is chain contractible ($\Z[X]$ denotes the free abelian group generated by $X$). In particular, it has trivial homology in each dimension. \end{lem}

The proof is standard (one constructs inductively a presimplicial contraction).

\begin{exmp} \label{magal1} Let $S \longrightarrow C$ be a $\P$-projective presimplicial resolution of $C$ and suppose $Q \in \P$. Then the augmented presimplicial set 
$$ \Hom_{\C}(Q,S) \longrightarrow \Hom_{\C}(Q,C)$$
satisfies the conditions of \ref{contract}. In particular, the homologies of the augmented chain complex
$$ K(\Z[\Hom_{\C}(Q, S)]) \longrightarrow \Z[\Hom_{\C}(Q,C)]$$
vanish.

\end{exmp}

Now suppose again that $\C$ is a category with finite limits, $\P$ a projective class in $\C$, $\B$ an abelian category, and $F \colon \C \longrightarrow \B$ a functor. Let $\Z\C$ denote the free preadditive category generated by $\C$ [M72, \S 1], i.e., the objects of $\Z\C$ are those of $\C$, and for any objects $C$ and $D$ in $\C$, $\Hom_{\Z\C}(C, D)$ is the free abelian group generated by $\Hom_{\C}(C, D)$. The composition of morphisms in $\Z\C$ is induced by that in $\C$. Clearly, $\C$ is a subcategory of $\Z\C$. Further, since the category $\B$ is abelian (and therefore additive), the functor $F \colon \C \longrightarrow \B$ can be uniquely extended to an additive functor
$$ F_{ad} \colon \Z\C \longrightarrow \B.$$
The following proposition relates the simplicial derived functors of $F$ to the Eilenberg-Moore derived functors of $F_{ad}$.

\begin{prop}\label{aditTV} Let $\C$ be a category with finite limits, $\P$ a projective class in $\C$, $\B$ an abelian category, and $F \colon \C \longrightarrow \B$ a functor. Then:

{\rm (i)} For any $\P$-projective presimplicial resolution $ S \longrightarrow C,$ the augmented chain complex
$$K(S) \longrightarrow C$$
in $\Z\C$ is a $\P$-resolution of $C$ in the sense of Definition \ref{res}.

{\rm (ii)} For any $C \in \C$, there is a natural isomorphism
$$\L^{\Delta}_nF(C) \cong \L_n^{\P}F_{ad}(C), \;\;\; n \geq 0. $$ \end{prop}

\begin{prf} The first claim immediately follows from \ref{magal1} and the definition of $\Z\C$. The second claim is a consequence of the first one and the definition of $F_{ad}$. Indeed, if $ S \longrightarrow C$ is a $\P$-projective presimplicial resolution of $C$, then we have
\begin{align*} \L^{\Delta}_nF(C) =H_n(K(F(S))) = \\
H_n(F_{ad}(K(S))) = L_n^{\P}F_{ad}(C). \end{align*} $\;$ \hfill \end{prf}

\begin{remk} Proposition \ref{aditTV} is essentially due to Barr and Beck [BB69, \S 5]. More precisely, in the case when the projective class $\P$ comes from a cotriple (see [TV69, \S 3]) the above statement is proved in [BB69, \S 5]. (The cotriple derived functor theory of Barr-Beck is a special case of the Tierney-Vogel theory [TV69 \S 3].) Thus \ref{aditTV} is a simple generalization of the result of Barr and Beck. \end{remk}

\vspace{0.2cm}

\section{Pseudocubical objects in idempotent complete preadditive categories}

\begin{defi}[{[P09, 2.2]}] A pseudocubical object X in a category $\C$ is a family of objects $(X_n \in \C)_{n \geq 0}$ together with face operators 
$$ \partial_i^0, \partial_i^1 \colon X_n \longrightarrow X_{n-1}, \;\;\; n \geq 1,\; \; 1 \leq i \leq n,$$
and pseudodegeneracy operators
$$ s_i \colon X_{n-1} \longrightarrow X_n, \;\;\;n \geq 1, \;\; 1 \leq i \leq n,$$
satisfying the pseudocubical identities 
\begin{align*}
\pa_i^\al\pa_j^\ve&=\pa_{j-1}^\ve\pa_i^\al\qquad
i<j,\quad\al,\ve\in\{0,1\},
\intertext{and}%
\pa_i^\al s_j&=\begin{cases}
s_{j-1}\pa_i^\al\ &i<j,\\
\id\ &i=j,\\
s_j\pa_{i-1}^\al\ &i>j,
\end{cases}
\end{align*}
for $\al\in\{0,1\}$.
\end{defi}

Important examples of pseudocubical objects appear in a natural way: Let $\C$ be a category with finite limits and $\P$ a projective class in $\C$. Then for any object $C \in \C$, there is a $\P$-exact augmented pseudocubical object
$$X \longrightarrow C$$
with $X_n \in \P$, $n \geq 0$, called $\P$-projective pseudocubical resolution of $C$ (see [P09, \S 3] for details).

In [P09] we use the normalized chain complex of a pseudocubical object in an abelian category to define the cubical derived functors. (Note that the normalized chain complex of a cubical object in an abelian category was introduced by \'Swiatek in [\'{S}75].) Below we recall the definition and some properties of the normalized chain complex of a pseudocubical object in the general setting of idempotent complete preadditive categories.
These are needed to prove a cubical analog of Proposition \ref{aditTV} in the next section.

\begin{defi} A preadditive category $\A$ is said to be idempotent complete if any idempotent $p :E \longrightarrow E$ in $\A$ (i.e., $p^2=p$) has a kernel. That is, there is a morphism
$$i \colon \Ker (p) \longrightarrow E$$
with $pi=0$, and for any morphism $f \colon F \lra E$, satisfying $pf=0$, there is a unique morphism $ g \colon F \lra \Ker (p)$ such that $ig=f$. \end{defi}

The following two propositions are well known (see e.g. [K78]).

\begin{prop} \label{split} Let $\A$ be an idempotent complete preadditive category and $p \colon E \lra E$ an idempotent in $\A$. Then there is a diagram
$$\xymatrix{ \Ker(p) \ar@<-0.5ex>[r]_-{i_1} & E \ar@<-0.5ex>[l]_-{\pi_1} \ar@<0.5ex>[r]^-{\pi_2} & \Ker(1-p) \ar@<0.5ex>[l]^-{i_2}}$$
such that 
\begin{align*}\pi_1i_1=1,\;\; \pi_2i_2=1,\\
\pi_1i_2=0,\;\; \pi_2i_1=0,\\
i_1\pi_1=1-p,\;\; i_2\pi_2=p. \end{align*}
In particular, the coproduct $\Ker(p) \oplus \Ker(1-p)$ exists in $\A$ and is isomorphic to $E$. \end{prop}

\begin{prop} Let $\A$ be a preadditive category. Then there exists an idempotent complete preadditive category $\wt{\A}$ and a full additive embedding
$$ \varphi \colon \A \lra \wt{\A}$$
satisfying the following universal property: For any idempotent complete preadditive category $\D$ and an additive functor $\psi : \A \lra \D$, there is an additive functor $\psi' \colon \wt{\A} \lra \D$ which makes the diagram 
$$\xymatrix{ \A \ar[rr]^-{\varphi} \ar[rd]_-{\psi} & & \wt{\A} \ar[ld]^-{\psi'} \\ & \D & }$$
commute up to a natural equivalence, and which is unique up to a natural isomorphism. \end{prop}

\;

Let $X$ be a pseudocubical object in an idempotent complete preadditive category $\D$.

\begin{defi} The unnormalized chain complex $C(X)$ associated to $X$ is defined by
$$ C(X)_n = X_n,  \;\;\; n \geq 0,$$
$$\partial= \sum_{i=1}^n(-1)^{i} (\partial_i^1-\pa_i^0) \colon C(X)_n \longrightarrow C(X)_{n-1}, \;\;\; n >0.$$ \end{defi}
The pseudocubical identities show that $\pa^2=0$. Moreover, they imply that the morphisms
$$\sigma_n^{X} = (1-s_1\pa_1^1)(1-s_2\pa_2^1) \cdots (1-s_n\pa_n^1) \colon X_n \lra X_n,\;\;\; n \geq 0,\;\;\; (\sigma_0=1)$$
are idempotents and form an endomorphism of the chain complex $C(X)$. We denote this endomorphism by
$$\sigma^{X} \colon C(X) \lra C(X).$$
Since $(\sigma^{X})^2 = \sigma^{X}$ and the category $\D$ is idempotent complete, the chain map $\sigma^{X}$ has a kernel $\Ker \sigma^{X}$ in the category of non-negative chain complexes in $\D$. Furthermore, by \ref{split}, there is a diagram in the category of chain complexes
$$\xymatrix{ \Ker(\sigma^{X}) \ar@<-0.5ex>[r]_-{i_1} & C(X) \ar@<-0.5ex>[l]_-{\pi_1} \ar@<0.5ex>[r]^-{\pi_2} & \Ker(1-\sigma^{X}) \ar@<0.5ex>[l]^-{i_2}}$$
such that
\begin{align*}\pi_1i_1=1,\;\; \pi_2i_2=1,\\
\pi_1i_2=0,\;\; \pi_2i_1=0,\\
i_1\pi_1=1-\sigma^{X},\;\; i_2\pi_2=\sigma^{X}. \end{align*}

\begin{defi} Let $X$ be a pseudocubical object in an idempotent complete preadditive category $\D$. The chain complex $\Ker(1-\sigma^{X})$, denoted by $N(X)$, is called the normalized chain complex of $X$. \end{defi}

\begin{remk} If $\D$ is an abelian category, then $N(X)$ admits the following description:

\begin{align*}
&N(X)_0=X_0,\quad N(X)_n=\us{i=1}{\os{n}{\cap}}\Ker(\pa_i^1),\quad n>0,\\
&\pa=\sum_{i=1}^n(-1)^{i+1}\pa_i^0:N(X)_n\to N(X)_{n-1},\quad n>0.
\end{align*}
Thus in the abelian case one does not need pseudodegeneracies to define $N(X)$. \end{remk}

Next, we recall the construction of cubical derived functors from [P09, \S 3]. Let $\C$ be a category with finite limits, $\P$ a projective class in $\C$, $\B$ an abelian category, and $F \colon \C \longrightarrow \B$ a functor. Then the cubical derived functors $\L^{\Box}_nF$ of $F$ with respect to the class $\P$ are defined as follows. For any object $C \in \C$, choose (once and for all) a $\P$-projective pseudocubical resolution 
$$ X \longrightarrow C$$
of $C$ and define 
$$\L^{\Box}_nF(C) = H_n(N(F(X))), \;\;\; n \geq 0.$$
\noindent The comparison theorem for precubical resolutions [P09, 3.3] and the homotopy invariance of the functor $N$ [P09, 3.6] imply that the objects $\L^{\Box}_nF(C)$ are well-defined and functorial in $F$ and $C$.

Note that one cannot use the unnormalized chain complex $C(X)$ instead of $N(X)$ to define the cubical derived functors [P09, 3.8].

The following lemma is the main technical tool for proving a cubical analog of Proposition \ref{aditTV}.

\begin{lem} \label{Nshenaxva} Let $F \colon \D \lra \D'$ be an additive functor between idempotent complete preadditive categories. Then for any pseudocubical object $X$ in $\D$, there is a natural isomorphism
$$F(N(X)) \cong N(F(X))$$
of chain complexes in $\D'$.
\end{lem}

\begin{prf} Applying the additive functor $F$ to the diagram
$$\xymatrix{ \Ker(\sigma^{X}) \ar@<-0.5ex>[r]_-{i_1} & C(X) \ar@<-0.5ex>[l]_-{\pi_1} \ar@<0.5ex>[r]^-{\pi_2} & \Ker(1-\sigma^{X})=N(X) \ar@<0.5ex>[l]^-{i_2},}$$
we get a digram in $\D'$
$$\xymatrix{ F(\Ker(\sigma^{X})) \ar@<-0.5ex>[r]_-{F(i_1)} & C(F(X)) \ar@<-0.5ex>[l]_-{F(\pi_1)} \ar@<0.5ex>[r]^-{F(\pi_2)} & F(N(X)) \ar@<0.5ex>[l]^-{F(i_2)}}$$
whose morphisms satisfy the following identities:

\begin{align*}F(\pi_1)F(i_1)=1,\;\; F(\pi_2)F(i_2)=1,\\
F(\pi_1)F(i_2)=0,\;\; F(\pi_2)F(i_1)=0,\\
F(i_1)F(\pi_1)=1-F(\sigma^{X}),\\
F(i_2)F(\pi_2)=F(\sigma^{X}). \end{align*}
Besides, it follows from the additivity of $F$ that $F(\sigma^{X})=\sigma^{F(X)}$, and hence we obtain
$$F(i_2)F(\pi_2)=\sigma^{F(X)}.$$
This finally implies that
$$F(N(X)) \cong \Ker(1-\sigma^{F(X)})=N(F(X)).$$ $\;$ \hfill \end{prf}

\vspace{0.2cm}

\section{Cubical derived functors and Eilenberg-Moore derived functors}

Let $\C$ be a category with finite limits, $\P$ a projective class, $\B$ an abelian category, and $F \colon \C \lra \B$ a functor. In this section we prove that for any object $C \in \C$, there is a natural isomorphism
$$\L^{\Box}_nF(C) \cong \L_n^{\P}F_{ad}(C), \;\;\; n \geq 0.$$
This together with \ref{aditTV} obviously implies Theorem \ref{mtavari}.

The proof of this isomorphism is similar to that of \ref{aditTV}. However, things become a little bit complicated in the cubical setting as we have to consider normalized chain complexes in order to get the ``right'' homology.

\begin{prop} Suppose $\A$ is a preadditive category, $\P$ a class of objects in $\A$, $\B$ an abelian category, and $F \colon \A \longrightarrow \B$ an additive functor. Suppose further that $\ol{\P}$ is the closure of the class $\P$ in the idempotent completion $\wt{\A}$, and $\wt{F} \colon \wt{\A} \lra \B$ the extension of $F$. Then for any $A \in \A$ which possesses a $\P$-resolution, there is a natural isomorphism
$$\L_n^{\P}F(A) \cong \L_n^{\overline{\P}}\wt{F}(A), \;\;\; n \geq 0.$$ \end{prop}

\begin{prf} Since $\wt{\A}$ has a zero object, any $\P$-resolution in $\A$ is a $\ol{\P}$-resolution in $\wt{\A}$. The rest follows from \ref{closure}. $\;$ \hfill \end{prf}

\begin{coro} \label{tilde} Assume that $\C$ ia a category with finite limits, $\P$ a projective class in $\C$, $\B$ an abelian category, and $F \colon \C \longrightarrow \B$ a functor. Assume further that $\wt{F_{ad}} \colon \wt{\Z\C} \lra \B$ is the extension of $F_{ad} \colon \Z\C \longrightarrow \B$ to the idempotent completion $\wt{\Z\C}$, and $\ol{\P}$ the closure of $\P$ in $\wt{\Z\C}$. Then for any object $C \in \C$, there is a natural isomorphism
$$\L_n^{\P}F_{ad}(C) \cong \L_n^{\overline{\P}}\wt{F_{ad}}(C), \;\;\; n \geq 0.$$ \end{coro}

Next we state the following technical

\begin{lem} \label{contractcube} Let $X \longrightarrow X_{-1}$ be an augmented pseudocubical set. Suppose that $\partial : X_0 \longrightarrow X_{-1}$ is surjective and the following conditions hold:

{\rm (i)} For any $x,y \in X_0$, satisfying $\partial x= \partial y$, there exists $z \in X_1$ such that $\pa_1^0z=x$ and $\pa_1^1z=y$.

{\rm (ii)} For any $n \geq 1$ and any collection of $2n +2$ elements $x_i^{\varepsilon} \in X_n$, $1 \leq i \leq n+1$, $\varepsilon \in \{0,1\}$, satisfying 
$$\pa_i^{\alpha}x_j^{\varepsilon} = \pa_{j-1}^{\varepsilon} x_i^{\alpha},\;\;\; 1 \leq i < j \leq n+1, \;\;\; \alpha, \varepsilon \in \{0,1\},$$
there exists $x \in X_{n+1}$, such that 
$$\pa_i^{\varepsilon}x=x_i^{\varepsilon},\;\;\;1 \leq i \leq n+1, \;\;\; \varepsilon \in \{0,1\}.$$
Then the augmented normalized chain complex
$$N(\Z[X]) \lra \Z[X_{-1}]$$
is chain contractible. In particular, it has trivial homology in each dimension. \end{lem}

We omit the routine details of the proof here. Note only that the main idea is to construct inductively a precubical homotopy equivalence between $X$ and the constant cubical object determined by $X_{-1}$ and then use the homotopy invariance of the functor $N$ [P09, 3.6].

\begin{exmp} \label{magal2} Let $X \longrightarrow C$ be a $\P$-projective pseudocubical resolution of $C$ and suppose $Q \in \P$. Then the augmented pseudocubical set 
$$ \Hom_{\C}(Q,X) \longrightarrow \Hom_{\C}(Q,C)$$
satisfies the conditions of \ref{contractcube}. In particular, the homologies of the augmented chain complex
$$ N(\Z[\Hom_{\C}(Q, X)]) \longrightarrow \Z[\Hom_{\C}(Q,C)]$$
vanish. \end{exmp}

We are now ready to prove the main result of this section.

\begin{prop} \label{aditcube} Let $\C$ be a category with finite limits, $\P$ a projective class in $\C$, $\B$ an abelian category, and $F \colon \C \longrightarrow \B$ a functor. Then:

{\rm (i)} For any $\P$-projective pseudocubical resolution $X \lra C$, the augmented chain complex 
$$N(X) \lra C$$
in the category $\wt{\Z\C}$ is a $\ol{\P}$-resolution of $C$ in the sense of \ref{res}. ($\ol{\P}$ is the closure of $\P$ in $\wt{\Z\C}$.)

{\rm (ii)} For any $C \in \C$, there is a natural isomorphism
$$\L^{\Box}_nF(C) \cong \L_n^{\P}F_{ad}(C), \;\;\; n \geq 0.$$ \end{prop}

\begin{prf} For all $n \geq 0$, $N(X)_n \in \ol{\P}$  since  $N(X)_n$ is a retract of $X_n$ and $\ol{\P}$ is closed under retracts. Further, by \ref{Nshenaxva}, one has a natural isomorphism of augmented chain complexes
$$\xymatrix{ \Hom_{\wt{\Z\C}}(Q, N(X)) \ar[r] \ar[d]^-{\cong} & \Hom_{\wt{\Z\C}}(Q, C) \ar[d]^{\idd} \\ N(\Z[\Hom_{\C}(Q, X)]) \ar[r] & \Z[\Hom_{\C}(Q,C)] }$$
for any $Q \in \P$. It follows from \ref{magal2} that the lower chain complex is acyclic and thus so is the upper one. Consequently, the augmented chain complex $N(X) \lra C$ in $\wt{\Z\C}$ is $\P$-acyclic or, equivalently, $\ol{\P}$-acyclic. This completes the proof of the first claim.

Let us prove the second claim. By \ref{tilde}, it suffices to get a natural isomorphism
$$\L^{\Box}_nF(C) \cong \L_n^{\ol{\P}}\wt{F_{ad}}(C).$$
Choose any $\P$-projective pseudocubical resolution $X \lra C$. The first claim together with \ref{Nshenaxva} gives
\begin{align*}\L_n^{\Box}F(C)= H_n(N(F(X))=H_n(N(\wt{F_{ad}}(X))) \cong \\
H_n(\wt{F_{ad}}(N(X))) = \L_n^{\ol{\P}}\wt{F_{ad}}(C). \end{align*} $\;$ \hfill \end{prf}

Clearly, Theorem \ref{mtavari} is an immediate consequence of \ref{aditTV} and \ref{aditcube}.

\vspace{0.2cm}

\section{Connection with topology}

In this section we briefly explain that Theorem \ref{mtavari} generalizes the well-known fact that the cubical and simplicial singular homologies of a topological space are naturally isomorphic. For the definition and basic properties of the cubical singular homology see [M80].

Let $\Top$ denote the category of topological spaces, and let $\Delta^n$, $n \geq 0$, be the standard $n$-simplex. The class $\P_{\Delta}$ of all possible disjoint unions of standard simplices is a projective class in $\Top$ in the sense of [TV69, \S 2]. (Moreover, in fact, it comes from a cotriple [BB69, (10.2)].) Indeed, for any space $Y$, the map
$$\bigsqcup_{\Delta^n \to Y, \atop  n \geq 0} \Delta^n \lra Y,$$
where the disjoint union is taken over all possible continuous maps $\Delta^n \lra Y$, $n \geq 0$, is a $\P_{\Delta}$-epimorphism. Consider the functor
$$F \colon \Top \lra \Ab,\;\;\; F(Y) = H_0^{\Delta}(Y,A)=\Z[\pi_0Y] \otimes A,$$
where $\Ab$ is the category of abelian groups, $H_*^{\Delta}(Y, A)$ the simplicial singular homology of $Y$ with coefficients in an abelian group $A$, and $\pi_0Y$ the set of path components of $Y$. It follows from [BB69, (10.2)] and [TV69, (3.1) Theorem] that there is a natural isomorphism
$$\L^{\Delta}_nF(Y) \cong H_n^{\Delta}(Y,A), \;\;\; n \geq 0,$$
where the simplicial derived functors are taken with respect to the projective class $\P_{\Delta}$ (cf. [R69], [R72]). We sketch the proof of this natural isomorphism along the lines of [BB69, (10.2)]. The standard cosimplicial object $\Delta^\bullet$ gives rise to an augmented simplicial functor
$$F_{\bullet} \lra F, \;\;\; F_n(Y)= \Z[\Hom_{\Top}(\Delta^n, Y)] \otimes A.$$
Further, suppose $S_{\bullet} \lra Y$ is a $\P_{\Delta}$-projective presimplicial resolution of $Y$. Evaluating $F_{\bullet}$ on $S_{\bullet}$ yields a bipresimplicial abelian group. It is easily seen that both resulting spectral sequences collapse at $E^2$. Finally, playing these two spectral sequences against each other gives the desired isomorphism.

Similarly, one can describe the cubical singular homologies $H_n^{\Box}(Y, A)$ as cubical derived functors of the functor $F(Y)=\Z[\pi_0Y] \otimes A$. For this one uses the class $\P_{\Box}$ consisting of all possible disjoint unions of standard cubes. The class $\P_{\Box}$ is a projective class in $\Top$ and there is a natural isomorphism
$$\L^{\Box}_nF(Y) \cong H_n^{\Box}(Y,A), \;\;\; n \geq 0,$$
where the cubical derived functors are taken with respect to $\P_{\Box}$. The proof of this isomorphism is technically a little bit complicated compared to its simplicial counterpart as one has to consider spectral sequences of bipseudocubical objects and take care of the normalizations.

Note that the class $\P= \P_{\Delta} \cup \P_{\Box}$ is also a projective class in $\Top$. Obviously, the simplicial derived functors with respect to the class $\P_{\Delta}$ are naturally isomorphic to the simplicial derived functors with respect to $\P$. On the other hand, the cubical derived functors with respect to the class $\P_{\Box}$ are naturally isomorphic to the cubical derived functors with respect to $\P$. Thus, by \ref{mtavari}, there is a natural isomorphism
$$\L^{\Delta}_nF(Y) \cong \L^{\Box}_nF(Y), \;\;\; n \geq 0,$$
for any topological space $Y$, i.e.,
$$H_n^{\Delta}(Y,A) \cong H_n^{\Box}(Y,A), \;\;\; n \geq 0.$$

\begin{tabular}{l}
{\sc Mathematisches Institut}\\
{\sc Universit\"at Bonn}\\
{\sc Endenicher Allee 60}\\
{\sc 53115 Bonn, Germany}\\
\\[-10pt]
\emph{E-mail address}: irpatchk@math.uni-bonn.de
\end{tabular}

\end{document}